\documentclass[12pt]{article}
\usepackage{amssymb, amsmath}

\newtheorem{thm}{Theorem}[section]

\newtheorem{proposition}{Proposition}[section]
\newtheorem{cor}{Corollary}[section]
\newcommand{\n}{\nonumber}

\newcommand{\vare}{\varepsilon}
\newcommand{\s}{\sigma}
\newcommand{\si}{\sigma_R }

\renewcommand{\o}{\omega}

\newcommand{\bb}{\begin{equation}}
\newcommand{\ee}{\end{equation}}
\newcommand{\bq}{\begin{eqnarray}}
\newcommand{\eq}{\end{eqnarray}}
\newcommand{\bqn}{\begin{eqnarray*}}
\newcommand{\eqn}{\end{eqnarray*}}

\begin{document}
\title{ Remarks on the  Liouville type problem in the stationary 3D Navier-Stokes equations}
 \author{Dongho Chae\\
\ \\
Department of Mathematics\\
Chung-Ang University\\
 Seoul 156-756, Republic of Korea\\
email: dchae@cau.ac.kr}
\date{}
\maketitle
\begin{abstract}
We study the Liouville type problem for the stationary 3D Navier-Stokes  equations on $\Bbb R^3$.
Specifically,  we prove that if $v$ is a smooth solution  to (NS) satisfying $\o={\rm curl}\,v \in L^q (\Bbb R^3) $ for some 
$\frac32 \leq q< 3$,  and $|v(x)|\to 0$
as $|x|\to +\infty$, then either $v=0$ on $\Bbb R^3$, or $\int_{\Bbb R^6} \Phi_+ dxdy=\int_{\Bbb R^6} \Phi_- dxdy=+\infty$, where
$\Phi(x,y) :=\frac{1}{4\pi}\frac{\o (x)\cdot(x-y)\times (v(y)\times \o(y) )}{|x-y|^3} $, and $\Phi_\pm:=\max\{ 0, \pm \Phi\}$. The proof uses crucially the structure of nonlinear term of the equations.
  \noindent{\bf
keywords:} stationary Navier-Stokes equations, Liouville type theorem
\end{abstract}

\section{Introduction}
 \setcounter{equation}{0}
 We consider the following stationary Navier-Stokes equations(NS) on $\Bbb R^3$.
 \bq\label{11}
 (v\cdot \nabla )v&=&-\nabla p+\Delta v,\\
 \label{12}
    \mathrm{div}\, v&=&0,
     \eq
  where $v(x)=(v_1(x), v_2(x), v_3(x))$ and $p=p(x)$ for all $ x\in \Bbb R^3$.
  The system is equipped with  the  boundary condition:
   \bb\label{13}
 |v(x)|\to 0 \quad \mbox{uniformly}\quad \mbox{as}\quad |x|\to +\infty.
 \ee  
 In addition to (\ref{13})  one usually  also assume following finite enstrophy condition.
 \bb\label{14}
 \int_{\Bbb R^3} |\nabla v|^2 dx <\infty,
 \ee
 which is physically natural.  It is well-known that any weak solution of (NS) satisfying (\ref{14}) is smooth. Actually, the regularity result for  the $L^\infty_t L^3_x$-weak solution of the non-stationary Navier-Stokes equations proved in \cite{es}  implies immediately that $v\in L^3(\Bbb R^3)$ is enough to guarantee the  regularity. A long standing open question  for  solution of (NS) satisfying the conditions (\ref{13}) and (\ref{14}) is that if it is trivial (namely, $v=0$ on $\Bbb R^3$), or not.  We refer the book by Galdi(\cite{gal})
 for the details on the motivations and historical backgrounds on the problem and the related results.  
 As a partial progress to the problem we mention that the condition $v\in L^{\frac92} (\Bbb R^3)$ implies that $v=0$ (see  Theorem X.9.5, pp. 729 \cite{gal}).  Another condition,  $\Delta v\in L^{\frac65}(\Bbb R^3)$ is also shown to imply $v=0$(\cite{cha}). 
 For studies on the Liouville type problem in the {\em non-stationary} Navier-Stokes equations, we refer \cite{ser}. Our aim in this paper is to prove the following:
  \begin{thm}
 Let $v$ be a smooth solution to  (NS) on $\Bbb R^3$ satisfying (\ref{13}). 
 Suppose there exists $q\in [\frac32, 3)$ such that  $\o \in  L^q(\Bbb R^3)$.
We set
\bb\label{int} \Phi (x,y):= \frac{1}{4\pi}\frac{\o (x)\cdot(x-y)\times (v(y)\times \o(y) )}{|x-y|^3} 
\ee
for all $(x,y)\in \Bbb R^3\times \Bbb R^3$ with $x\not=y$, and define
 $$\Phi_+ (x,y):=\max\{ 0, \Phi (x,y)\}, \quad \Phi_{-}(x,y):=\max\{ 0, -\Phi(x,y)\}.
 $$
 Then, either 
 \bb\label{zer}
 v=0 \quad\mbox{on}\quad \Bbb R^3,
 \ee
 or
 \bb\label{inf}
\int_{\Bbb R^3}\int_{\Bbb R^3} \Phi_+ (x,y) dxdy =\int_{\Bbb R^3}\int_{\Bbb R^3} \Phi_- (x,y) dxdy =+\infty.
\ee
\end{thm} 
 \noindent{\em Remark 1.1 }  One can show that if  $\o\in L^{\frac95} (\Bbb R^3)$  is satisfied together with (\ref{13}), then  (\ref{zer}) holds.
In order to see this  we first recall the estimate of the Riesz potential on  $\Bbb R^3$(\cite{ste}),
\bb\label{rie} \| I_\alpha (f)\|_{L^q}\leq C\|f\|_{L^p} , \quad \frac{1}{q}=\frac{1}{p}-\frac{\alpha}{3},\quad 1\leq p<q<\infty,
\ee
where
 $$I_\alpha (f):=C\int_{\Bbb R^3} \frac{f(y)}{|x-y|^{3-\alpha}} dy, \quad 0<\alpha <3 $$
 for a positive constant $C=C(\alpha)$.
Applying (\ref{rie}) with $\alpha=1$, we obtain by   the H\"{o}lder inequality,
 \bqn
\lefteqn{ \int_{\Bbb R^3} \int_{\Bbb R^3} |\Phi (x,y)| dydx \leq \int_{\Bbb R^3} \int_{\Bbb R^3} \frac{ |\o (x)||\o(y)||v(y)|}{|x-y|^2} dydx }\n \\
&&\quad \leq \left(\int_{\Bbb R^3} |\o(x)|^{\frac95} dx \right)^{\frac59}\left\{\int_{\Bbb R^3}\left( \int_{\Bbb R^3} \frac{|\o(y)||v(y)|}{|x-y|^2} dy\right)^{\frac94} dx \right\}^{\frac49}\n \\
 && \quad\leq C\|\o\|_{L^{\frac95}} \left( \int_{\Bbb R^3} |\o|^{\frac97} |v|^{\frac97}dx\right)^{\frac79}\n \\
 &&\quad \leq C\|\o\|_{L^{\frac95}} \left( \int_{\Bbb R^3} |\o|^{\frac95} dx\right) ^{\frac59} \left(\int_{\Bbb R^3} |v|^{\frac92}dx\right)^{\frac29} \n \\
&& \quad\leq C\|\o\|_{L^{\frac95}}^2 \|\nabla v\|_{L^{\frac95}} \leq C\|\o\|_{L^{\frac95}} ^3<+\infty,
 \eqn
 where we used  the Sobolev and the Calderon-Zygmund inequalities  
 \bb\label{so}
 \|v\|_{L^{\frac92}} \leq C \|\nabla v\|_{L^{\frac95}} \leq C \|\o\|_{L^{\frac95}}
 \ee 
 in the last step.   Thus, by the Fubini-Tonelli theorem, (\ref{inf}) cannot hold, and we are lead to (\ref{zer}) by application of the above theorem.   We  note that by (\ref{so})
 the condition $\o\in L^{\frac95} (\Bbb R^3)$, on the other hand, implies the previously known sufficient condition 
 $v\in L^{\frac92} (\Bbb R^3)$ of \cite{gal} mentioned above.

  \section{Proof of the main theorem}
 \setcounter{equation}{0}
 We first establish  integrability conditions on the vector fields for the Biot-Savart's  formula  in $\Bbb R^3$.
\begin{proposition}
 Let $\xi=\xi(x)=(\xi_1(x), \xi_2 (x), \xi_3 (x))$ and $\eta=\eta(x)=(\eta_1 (x), \eta_2 (x), \eta_3 (x)) $ be smooth  vector fields on $\Bbb R^3$.
 Suppose there exists $q \in [1, 3)$ such that $\eta\in L^q(\Bbb R^3)$.   Let $\xi$ solve
\bb\label{pro1}
\Delta \xi=-\nabla \times \eta,
\ee
under the boundary condition; either
 \bb\label{pro2}
 |\xi (x)|\to 0 \quad \mbox{uniformly}\quad \mbox{as}\quad |x|\to +\infty, 
\ee
or 
\bb\label{pro3}
 \xi \in L^s(\Bbb R^3) \quad \mbox{for some}\quad s\in [1,\infty).
\ee
Then, the solution of (\ref{pro1})  is given by 
\bb\label{21}
\xi(x)=\frac{1}{4\pi}\int_{\Bbb R^3} \frac{(x-y)\times \eta(y)}{|x-y|^3}dy \quad \forall x\in \Bbb R^3.
 \ee
 \end{proposition}
\noindent{\bf Proof   } We introduce a cut-off function $\sigma\in C_0
^\infty(\Bbb R^N)$ such that
$$
   \sigma(|x|)=\left\{ \aligned
                  &1 \quad\mbox{if $|x|<1$},\\
                     &0 \quad\mbox{if $|x|>2$},
                      \endaligned \right.
$$
and $0\leq \sigma (x)\leq 1$ for $1<|x|<2$.  For each $R>0$ we
define $\s_R (x):= \s \left(\frac{|x|}{R}\right)$.  Given $\vare >0$ we denote 
$B_\vare (y)=\{ x\in \Bbb R^3\, |\, |x-y|<\vare\}$. 
Let us fix $y\in \Bbb R^3$ and $\vare \in (0, \frac{R}{2})$.
We multiply  (\ref{pro1}) by $\frac{\s_R (|x-y|)}{|x-y|},$ 
and integrate it with respect to the variable $x$ over 
$\Bbb R^3\setminus B_{\vare} (y)$. Then, 
\bb\label{23}
 \int_{\{|x-y|>\vare\}}  \frac{\Delta \xi\,\si}{|x-y|} dx =  -\int_{\{|x-y|>\vare\}}  \frac{\si \nabla \times \eta (y)}{|x-y|} dx .
 \ee
Since $\Delta \frac{1}{|x-y|}=0$ on $\Bbb R^3\setminus B_\vare (y)$, one has
\bqn
\lefteqn{ \frac{\Delta \xi\si}{|x-y|}= \sum_{i=1}^3\partial_{x_i} \left( \frac{\partial_{x_i} \xi\si}{|x-y|} \right) -\sum_{i=1}^3\partial_{x_i}\left( \frac{\xi \partial_{x_i}\si}{|x-y|}\right)}
 \n \\
&& -\sum_{i=1}^3\partial_{x_i} \left( \xi \si \partial_{x_i} (\frac{1}{|x-y|} )\right)
+\frac{\xi\Delta \si}{|x-y|} +2\sum_{i=1}^3\xi \partial_{x_i}(\frac{1}{|x-y|}) \partial_{x_i}\si .
 \eqn
 Therefore, applying the divergence theorem,  and observing $\partial_{\nu} \si=0$ on $\partial B_\vare (y)$, we have
 \bq\label{24}
\lefteqn{ \int_{\{|x-y|>\vare\}}  \frac{\Delta \xi\si}{|x-y|} dx = \int_{\{|x-y|=\vare\}} \frac{ \partial_\nu \xi}{|x-y|}   dS}
\n \\ 
&&\qquad -\int_{\{ |x-y|=\vare\}} \frac{\xi}{|x-y|^2} dS + \int_{\{ |x-y|>\vare\}} \frac{\xi \Delta \si}{|x-y|} dx \n \\
&&\qquad-2\int_{\{ |x-y|>\vare\}} \frac{ (x-y)\cdot \nabla \si \,\xi}{|x-y|^3} dx
  \eq
  where $\partial_\nu(\cdot)$ denotes  the outward normal derivative on $\partial B_\vare(y)$. Passing $\vare \to 0$, one  can  easily compute 
 that 
\bq\label{25}
\mbox{\rm RHS of (\ref{24}) }&\to& -4\pi \xi (y) +\int_{\Bbb R^3 } \frac{\xi\Delta \si}{|x-y|} dx  -2\int_{\Bbb R^3} \frac{ (x-y)\cdot \nabla \si \,\xi}{|x-y|^3} dx \n\\
&:=& I_1+I_2+I_3.
\eq
Next, using the formula
  $$
   \frac{\si \nabla \times \eta}{|x-y|} =\nabla \times \left( \frac{\si \eta}{|x-y|}\right) -\frac{\nabla \si \times \eta}{|x-y|} +
   \frac{(x-y)\times \eta\si }{|x-y|^3} ,
 $$
   and using the divergence theorem, we obtain the following representation for the right hand side of (\ref{23}).
   \bq\label{26}
   \lefteqn{ \int_{\{|x-y|>\vare\}}  \frac{\si \nabla \times \eta}{|x-y|}  dx  =\int_{\{|x-y|=\vare\}}  \nu \times \left( \frac{\eta}{|x-y|}\right) dS
   }\n \\
   && - \int_{\{|x-y|>\vare\}}\frac{\nabla \si \times \eta}{|x-y|} dx+ \int_{\{|x-y|>\vare\}} \frac{(x-y)\times \eta\si }{|x-y|^3}dx, \n \\
   \eq
   where we denoted $\nu=\frac{y-x}{|y-x|}$,  the outward unit normal vector on $\partial B_\vare (y)$. Passing $\vare\to 0$, we easily deduce
 \bq \label{27}
 \mbox{\rm RHS of (\ref{26}) }&\to &- \int_{\Bbb R^3}\frac{\nabla \si \times \eta}{|x-y|} dx+ \int_{\Bbb R^3} \frac{(x-y)\times \eta\si }{|x-y|^3}dx\n \\
 &:=&J_1+J_2 \quad \mbox{as} \quad \vare \to 0.
 \eq   
 We now pass $R\to \infty$ for each term of (\ref{25}) and (\ref{27}) respectively below.
 Under the boundary condition (\ref{pro2}) we estimate:
 \bqn
 |I_2|&\leq& \int_{\{ R\leq |x-y|\leq 2R\}}  \frac{|\xi(x)| |\Delta\si (x-y)|}{|x-y|} dx\n \\
 &\leq& \frac{\|\Delta \s\|_{L^\infty}}{R^2}   \sup_{R\leq |x|\leq 2R} |\xi(x)| \left(\int_{\{R\leq |x-y|\leq 2R\}} dx \right)^
 {\frac{2}{3}}  \left(\int_{\{R\leq |x-y|\leq 2R\}} \frac{dx}{ |x-y|^3} \right)^
 {\frac{1}{3}} \n \\
 &\leq& C \|\Delta \s\|_{L^\infty} \left(\int_{R} ^{2R} \frac{dr}{r} \right) ^{\frac23}\sup_{R\leq |x-y|\leq 2R} |\xi(x)| \to 0
 \eqn 
 as $R\to \infty$ by the assumption (\ref{pro2}), while under the condition (\ref{pro3}) we have
 \bqn
 |I_2|&\leq& \int_{\{ R\leq |x-y|\leq 2R\}}  \frac{|\xi(x)| |\Delta\si (x-y)|}{|x-y|} dx\n \\
 &\leq& \frac{\|\Delta \s\|_{L^\infty}}{R^2}  \|\xi\|_{L^s }\left(\int_{\{0\leq |x-y|\leq 2R\}} \frac{dx}{ |x-y|^{\frac{s}{s-1}}}  \right)^
 {\frac{s-1}{s}} \n \\
 &\leq& C R^{-\frac{3}{s}} \|\Delta \s\|_{L^\infty}   \|\xi\|_{L^s }\to 0
 \eqn  
 as $R\to \infty$. Similarly, under (\ref{pro2}) 
  \bqn
 |I_3|&\leq& 2 \int_{\{ R\leq |x-y|\leq 2R\}}  \frac{|\xi(x)| |\nabla\si (x-y)|}{|x-y|^2} dx\n \\ 
 &\leq& \frac{C\|\nabla \s\|_{L^\infty}}{R}  
 \sup_{R\leq |x|\leq 2R} |\xi(x)| \left(\int_{\{R\leq |x-y|\leq 2R\}} dx \right)^
 {\frac{1}{3}}  \left(\int_{\{R\leq |x-y|\leq 2R\}} \frac{dx}{ |x-y|^3}  \right)^
 {\frac{2}{3}}   \n \\ 
 &\leq & C\|\nabla \s\|_{L^\infty}\left(\int_{R} ^{2R} \frac{dr}{r} \right) ^{\frac23} \sup_{R\leq |x-y|\leq 2R} |\xi(x)|\to 0 
 \eqn 
 as $R\to \infty$, while under the condition (\ref{pro3}) we estimate
   \bqn
 |I_3|&\leq& 2 \int_{\{ R\leq |x-y|\leq 2R\}}  \frac{|\xi(x)| |\nabla\si (x-y)|}{|x-y|^2} dx\n \\ 
 &\leq& \frac{C\|\nabla \s\|_{L^\infty}}{R}  
  \|\xi\|_{L^s(R\leq |x-y|\leq 2R) } \left(\int_{\{0\leq |x-y|\leq 2R\}} \frac{dx}{ |x-y|^{\frac{2s}{s-1}}} \right)^
 {\frac{s-1}{s}}   \n \\ 
 &\leq & C R^{-\frac{3}{s}} \|\nabla \s\|_{L^\infty} \|\xi\|_{L^s }\to 0 
 \eqn  
 as $R\to \infty$.
 Therefore, 
the right hand side  of (\ref{24}) converges to $-4\pi \xi(y)$  as $R\to \infty.$  For  $J_1, J_2$ we estimate
 \bqn
 |J_1|&\leq& \int_{\{ R\leq |x-y|\leq 2R\}}  \frac{|\nabla \si| |\eta|}{|x-y|} dx\n \\
 &\leq & \frac{C\|\nabla \s\|_{L^\infty}}{R}  
 \|\eta\|_{L^{q}(R\leq |x-y|\leq 2R)}\left(\int_{\{0\leq |x-y|\leq 2R\}}  \frac{dx}{|x-y|^{\frac{q}{q-1}}} \right)^
 {\frac{q-1}{q}}\n \\
 &\leq& C\|\nabla \s\|_{L^\infty}  \|\eta\|_{L^{q} (R\leq |x-y|\leq 2R)}  R^{-\frac{2}{q}}\to 0
 \eqn 
 as $R\to \infty$.  
  In pasing $R\to \infty$ in $J_2$ of (\ref{27}), in order to use the dominated convergence theorem, we estimate
\bq\label{j1}
 \int_{\Bbb R^3} \left|\frac{(x-y)\times \eta (y)}{|x-y|^3}\right|dx 
 &\leq&  \int_{\{|x-y|<1\}}\frac{|\eta| }{|x-y|^2}dx
 +  \int_{\{|x-y| \geq 1\}} \frac{|\eta| }{|x-y|^2}dx \n \\
 &:=&J_{21}+J_{22}.
 \eq
$J_{21}$ is easy to handle as follows.
\bb\label{j2}
 J_{21} \leq  \|\eta \|_{L^\infty (B_1(y))}\int_{\{ |x-y|<1\}} \frac{dx}{|x-y|^2} 
 = 4\pi\|\eta\|_{L^\infty (B_1(y))}<+\infty.
\ee
   For $J_{22}$ we estimate
 \bq\label{j3}
 J_{22}&\leq&  \left(\int_{\Bbb R^3} |\eta|^q dx \right)^{\frac{1}{q} }  \left( \int_{ \{ |x-y|>1\}} \frac{dx}{|x-y|^{\frac{2q}{q-1}} }\right)^{\frac{q-1}{q}} \n \\
 &\leq& C \|\eta \|_{L^{q}}\left( \int _{1} ^\infty r ^{\frac{-2}{q-1}} dr \right) ^{\frac{q-1}{q}}<+\infty,
 \eq
 if $1<q<3$.  In the case of $q=1$ we estimate simply
\bb\label{j4} J_{22} \leq  \int_{ \{ |x-y|>1\}} |\eta| dx \leq \|\eta\|_{L^1}.
 \ee
 Estimates of (\ref{j1})-(\ref{j4}) imply 
 $$ \int_{\Bbb R^3} \left|\frac{(x-y)\times \eta (y)}{|x-y|^3}\right|dx <+\infty.
 $$
  Summarising the above computations, one can pass first $\vare \to 0$, and then $R\to +\infty$ in 
 (\ref{23}), applying the dominated convergence theorem, to obtain finally (\ref{21}). $\square$\\
  \ \\
\begin{cor}
Let $v$ be a smooth solution to (\ref{11})-(\ref{13}) such that $\o \in L^q(\Bbb R^3)$ for some $q\in [\frac32, 3)$.
Then, we have
\bb\label{bio}
v(x)=\frac{1}{4\pi} \int_{\Bbb R^3} \frac{(x-y)\times \o(y)}{|x-y|^3} dy,
\ee
and
\bb\label{vor}
\o (x)=\frac{1}{4\pi} \int_{\Bbb R^3} \frac{(x-y)\times (v(y)\times \o(y))}{|x-y|^3} dy.
\ee
\end{cor}
{\bf Proof }  Taking curl of the defining equation of the vorticity, $\nabla \times v=\o$, using div $v=0$, we have
$$ \Delta v=-\nabla \times \o, $$
which provides us with (\ref{bio}) immediately by application of Proposition 2.1. In order to show (\ref{vor}) we recall that,
using the vector identity $\frac12 \nabla |v|^2 =(v\cdot \nabla )v +v\times (\nabla \times v)$, one can rewrite (\ref{11})-(\ref{12}) as
 $$-v\times \o =-\nabla\left(p+\frac12 |v|^2\right) +\Delta v.$$ 
Taking curl on this, we obtain
$$
 \Delta \o= -\nabla \times (v\times \o).
$$ 
The formula (\ref{vor}) is deduced immediately  from this equations by applying the proposition 2.1.
For the allowed rage of $q$ we recall  the Sobolev and the Calderon-Zygmund inequalities(\cite{ste}),
 \bb\label{sob}
  \|v\|_{L^{\frac{3q}{3-q}}} \leq C\|\nabla v\|_{L^q}\leq C \|\o\|_{L^q}, \quad 1<q<3, 
\ee
which imply $v\times \o \in L^{\frac{3q}{6-q}} (\Bbb R^3)$ if $\o \in L^q(\Bbb R^3)$. 
We also note that  $\frac32 \leq q<3$ if and only if $1\leq \frac{3q}{6-q} <3$. $\square$ \\
\ \\
 \noindent{\bf Proof of Theorem 1.1 }  
 Under the hypothesis (\ref{13}) and $\o \in L^q(\Bbb R^3)$ with $q\in [\frac32, 3)$ both of the relations  (\ref{bio}) and (\ref{vor}) are valid.
We first prove the following.\\

\noindent{\em Claim: } For each $x, y\in \Bbb R^3$
\bb\label{th1a}
0\leq |\o(x)|^2=\int_{\Bbb R^3} \Phi (x,y)dy \leq \int_{\Bbb R^3} |\Phi (x,y)| dy <+\infty,
\ee
and
\bb\label{th1b}
0=\int_{\Bbb R^3} \Phi (x,y)dx \leq \int_{\Bbb R^3} |\Phi (x,y)| dx<+\infty.
\ee
{\em Proof of the claim: }  We  verify the following:
\bb\label{veri}
\int_{\Bbb R^3} |\Phi(x,y)|dy +\int_{\Bbb R^3}|\Phi (x,y)|dx <\infty \quad \forall (x,y)\in \Bbb R^3\times \Bbb R^3.
\ee
Decomposing the integral, and using the H\"{o}older inequality, we estimate
\bq\label{pf0}
 \int_{\Bbb R^3} |\Phi (x,y)|dy&\leq& |\o(x)|\left( \int_{\{ |x-y|\leq 1\}}  \frac{|v(y)||\o(y)|}{|x-y|^2} dy
                            +\int_{\{ |x-y|>1\}} \frac{|v(y)||\o(y)|}{|x-y|^2} dy \right)\n \\
                         &\leq& |\o(x)|\|v\|_{L^\infty(B_1 (x)) }\|\o\|_{L^\infty(B_1 (x))} \int_{\{ |x-y|\leq 1\}}  \frac{dy}{|x-y|^2}  \n \\
                            &&\qquad + |\o(x)|\|v\|_{L^{\frac{3q}{3-q}}} \|\o\|_{L^q} \left(\int_{\{ |x-y|\geq 1\}}  \frac{dy}{|x-y|^{\frac{6q}{4q-6}}}\right)^{\frac{4q-6}{3q}} \n \\
   & \leq &C |\o(x)|\|v\|_{L^\infty(B_1 (x)) }\|\o\|_{L^\infty(B_1 (x))}   \n \\
   &&\qquad + C|\o(x)| \|\o\|_{L^q}^2 \left(\int_1 ^\infty r^{\frac{q-6}{2q-3}} dr\right)^{\frac{4q-6}{3q}} <+\infty,                     
   \eq
 where we used (\ref{sob}) and the fact that $\frac{q-6}{3q-3} <-1$ if $\frac32 < q<3$. In the case $q=\frac32$ we estimate, instead,
 \bq\label{pf01}
 \int_{\Bbb R^3} |\Phi (x,y)|dy&\leq& |\o(x)|\left( \int_{\{ |x-y|\leq 1\}}  \frac{|v(y)||\o(y)|}{|x-y|^2} dy
                            +\int_{\{ |x-y|>1\}} \frac{|v(y)||\o(y)|}{|x-y|^2} dy \right)\n \\
                            &\leq& |\o(x)|\|v\|_{L^\infty(B_1 (x)) }\|\o\|_{L^\infty(B_1 (x))}  + |\o(x)|\|v\|_{L^{3}} \|\o\|_{L^{\frac32}} <+\infty.\n \\
 \eq
 We also have
 \bq\label{pf0a}
\lefteqn{ \int_{\Bbb R^3} |\Phi (x,y)|dx\leq |v(y)| |\o(y)|\left( \int_{\{ |x-y|\leq 1\}}  \frac{|\o(x)|}{|x-y|^2} dx
                            +\int_{\{ |x-y|>1\}} \frac{|\o(x)|}{|x-y|^2} dx\right)}\n \\
 &&\leq C |v(y)| |\o(y)| \| \o\|_{L^\infty (B_1 (y))}+|v(y)| |\o(y)| \|\o\|_{L^q} \left(\int_{\{ |x-y|>1\}} \frac{dx}{|x-y|^{\frac{2q}{q-1}} }\right)^{\frac{q-1}{q}} \n \\
 &&\leq C |v(y)| |\o(y)| \| \o\|_{L^\infty (B_1 (y))}+C|v(y)| |\o(y)| \|\o\|_{L^q}   \left(\int_1 ^\infty r^{-\frac{2}{q-1}} dr \right) ^{\frac{q-1}{q}}<+\infty,\n \\
 \eq
 where we used the fact that $ -\frac{2}{q-1} <-1$ if $ \frac32\leq q<3$.
  From (\ref{vor}) we immediately obtain
 \bq\label{pf1}
 \int_{\Bbb R^3} \Phi (x,y)dy &=& \o (x)\cdot \left(\frac{1}{4\pi}\int_{\Bbb R^3} \frac{(x-y)\times (v(y)\times \o(y))}{|x-y|^3} dy \right)\n \\
 &=& |\o (x)|^2 \geq 0, \quad \forall x\in \Bbb R^3
 \eq
and combining this with (\ref{pf0}),  we deduce (\ref{th1a}).  On the other hand, using (\ref{bio}), we find
\bq\label{pf2}
 \int_{\Bbb R^3} \Phi (x,y)dx&=& \frac{1}{4\pi} \int_{\Bbb R^3} \frac{\o (x)\cdot (x-y)\times (v(y)\times \o(y))}{|x-y|^3}dx \n \\
 &=& \left(\frac{1}{4\pi} \int_{\Bbb R^3} \frac{\o (x)\times(x-y)}{|x-y|^3}dx\right) \cdot v(y)\times \o(y)\n \\
 &=& v(y)\cdot  v(y)\times \o(y) =0
  \eq
 for all $y\in \Bbb R^3$, and combining this with (\ref{pf0a}),  we have proved  (\ref{th1b}).
 This completes the proof of the claim.\\
 \ \\
By the Fubini-Tonelli theorem we have
\bb\label{pos}
\int_{\Bbb R^3} \int_{\Bbb R^3} \Phi_+ (x,y) dxdy
=\int_{\Bbb R^3} \int_{\Bbb R^3} \Phi_+ (x,y) dydx:=\mathcal{I}_+ , 
\ee
and 
\bb
\label{neg}
\int_{\Bbb R^3} \int_{\Bbb R^3} \Phi_- (x,y) dxdy
=\int_{\Bbb R^3} \int_{\Bbb R^3} \Phi_- (x,y) dydx :=\mathcal{I}_-.
\ee
If (\ref{inf}) does not hold, then at least one of the two integrals $\mathcal{I}_+ , \mathcal{I}_-$ is finite. In this case, using (\ref{pos}) and (\ref{neg}),  we  can interchange the order of integrations in repeated integral as follows.
\bq\label{com}
\int_{\Bbb R^3}\int_{\Bbb R^3} \Phi (x,y)dxdy &=&\int_{\Bbb R^3} \int_{\Bbb R^3} \Phi_+ (x,y) dxdy
-\int_{\Bbb R^3} \int_{\Bbb R^3} \Phi_- (x,y) dxdy\n \\
&= &\int_{\Bbb R^3} \int_{\Bbb R^3} \Phi_+ (x,y) dydx
-\int_{\Bbb R^3} \int_{\Bbb R^3} \Phi_- (x,y) dydx\n \\
& =&\int_{\Bbb R^3}\int_{\Bbb R^3} \Phi (x,y)dydx.
\eq
Therefore,  from (\ref{pf1}) and (\ref{pf2}) combined with (\ref{com}) provide us with 
\bqn
\int_{\Bbb R^3} |\o(x)|^2 dx&=&\int_{\Bbb R^2} \int_{\Bbb R^3} \Phi (x,y) dydx
=\int_{\Bbb R^2} \int_{\Bbb R^3} \Phi (x,y) dxdy
=0.
\eqn
Hence, 
\bb\label{last}
\o=0\quad\mbox{on}\quad \Bbb R^3.
\ee
 We remark parenthetically  that in deriving (\ref{last}) it is not necessary to assume that $\int_{\Bbb R^3} |\o(x)|^2 dx <+\infty $, and we do not need to restrict ourselves to  $\o \in L^2 (\Bbb R^3)$.
Hence,   from (\ref{bio}) and (\ref{last}), we we conclude $v=0$ on $\Bbb R^3$. 
$\square$

      $$\mbox{\bf Acknowledgements}$$
This research was partially supported by NRF grants 2006-0093854 and  2009-0083521. 
  
\end{document}